\documentclass[a4paper,twoside,12pt]{amsart}
\RequirePackage{amsmath,amssymb,xypic}
\theoremstyle{plain}
\newtheorem{theorem}{Theorem}[section]
\newtheorem{corollary}[theorem]{Corollary}
\newtheorem{proposition}[theorem]{Proposition}
\newtheorem{lemma}[theorem]{Lemma}

\theoremstyle{definition}
\newtheorem{definition}[theorem]{Definition}
\newtheorem{example}[theorem]{Example}
\newtheorem{remark}[theorem]{Remark}

\numberwithin{equation}{section}

\newcommand{\Union}{\bigcup\limits}

\newcommand{\C}{\mathbb{C}}
\newcommand{\N}{\mathbb{N}}

\newcommand{\Z}{\mathbb{Z}}

\newcommand{\op}{\mathrm{op}}

\DeclareMathOperator{\id}{id}

\renewcommand{\to}[1][]{\xrightarrow[#1]{}}
\newcommand{\from}[1][]{\xleftarrow[#1]{}}

\newcommand{\Endo}[1][]{\mathrm{End}_{\raise1.5ex\hbox to.1em{}#1}}

\newcommand{\Hom}[1][]{\mathrm{Hom}_{\raise1.5ex\hbox to.1em{}#1}}

\newcommand{\RHom}[1][]{\mathrm{RHom}_{\raise1.5ex\hbox to.1em{}#1}}

\newcommand{\Ext}[2][]{\mathrm{Ext}_{\raise1.5ex\hbox to.1em{}#1}^{#2}}

\newcommand{\THom}[1][]{\mathrm{THom}_{\raise1.5ex\hbox to.1em{}#1}}

\newcommand{\Tens}[1][]{\mathbin{\otimes_{\raise1.5ex\hbox to-.1em{}#1}}}

\newcommand{\LTens}[1][]{\mathbin{\otimes_{\raise1.5ex\hbox to-.1em{}#1}^{L}}}

\newcommand{\Tor}[2][]{\mathrm{Tor}^{\raise1.5ex\hbox to.1em{}#1}_{#2}}

\def\sha{\mathcal{A}}
\def\shb{\mathcal{B}}

\def\she{\mathcal{E}}

\def\shl{\mathcal{L}}
\def\shm{\mathcal{M}}

\def\shp{\mathcal{P}}
\def\shq{\mathcal{Q}}
\def\shr{\mathcal{R}}

\def\shw{\mathcal{W}}

\newcommand{\shendo}[1][]{{\mathcal{E}nd}_{\raise1.5ex\hbox to.1em{}#1}}

\renewcommand{\hom}[1][]{{\mathcal{H}om}_{\raise1.5ex\hbox to.1em{}#1}}

\newcommand{\rhom}[1][]{{R\mathcal{H}om}_{\raise1.5ex\hbox to.1em{}#1}}

\newcommand{\ext}[2][]{{\mathcal{E}xt}_{\raise1.5ex\hbox to.1em{}#1}^{#2}}

\newcommand{\thom}[1][]{{T\mathcal{H}om}_{\raise1.5ex\hbox to.1em{}#1}}

\newcommand{\tens}[1][]{\mathbin{\otimes_{\raise1.5ex\hbox to-.1em{}#1}}}

\newcommand{\ltens}[1][]{\mathbin{\otimes_{\raise1.5ex\hbox to-.1em{}#1}^{L}}}

\newcommand{\tor}[2][]{{\mathcal{T}or}^{\raise1.5ex\hbox to.1em{}#1}_{#2}}

\newcommand{\oim}[1]{{#1}_*}

\newcommand{\opb}[1]{#1^{-1}}

\newcommand{\GHom}[1][]{\mathrm{GHom}_{\raise1.5ex\hbox to.1em{}#1}}

\newcommand{\GExt}[2][]{\mathrm{GExt}_{\raise1.5ex\hbox to.1em{}#1}^{#2}}

\newcommand{\FHom}[1][]{\mathrm{FHom}_{\raise1.5ex\hbox to.1em{}#1}}

\newcommand{\ghom}[1][]{{\mathcal{GH}om}_{\raise1.5ex\hbox to.1em{}#1}}

\newcommand{\gext}[2][]{{\mathcal{GE}xt}_{\raise1.5ex\hbox to.1em{}#1}^{#2}}

\newcommand{\fhom}[1][]{{\mathcal{FH}om}_{\raise1.5ex\hbox to.1em{}#1}}

\newcommand{\tenstop}[1][]{\mathbin{\hat{\otimes}_{\raise1.5ex\hbox to-.1em{}#1}}}

\newcommand{\homtop}[1][]{\mathcal{L}_{\raise1.5ex\hbox to.1em{}#1}}

\newcommand{\Homtop}[1][]{\mathrm{L}_{\raise1.5ex\hbox to.1em{}#1}}

\renewcommand{\O}{\mathcal{O}}

\def\absdoim#1{\underline{#1}_*}
\def\reldoim[#1]#2{\underline{#2}_{|{#1}*}}
\def\doim{\@ifnextchar [{\reldoim}{\absdoim}}

\def\absdeim#1{\underline{#1}_*}
\def\reldeim[#1]#2{\underline{#2}_{|{#1}*}}
\def\deim{\@ifnextchar [{\reldeim}{\absdeim}}

\def\absdopb#1{\underline{#1}^{-1}}
\def\reldopb[#1]#2{\underline{#2}_{|{#1}}^{-1}}
\def\dopb{\@ifnextchar [{\reldopb}{\absdopb}}

\def\absboim#1{\underline{\underline{#1}}_*}
\def\relboim[#1]#2{\underline{\underline{#2}}_{|{#1}*}}
\def\boim{\@ifnextchar [{\relboim}{\absboim}}

\def\absbeim#1{\underline{\underline{#1}}_*}
\def\relbeim[#1]#2{\underline{\underline{#2}}_{|{#1}*}}
\def\beim{\@ifnextchar [{\relbeim}{\absbeim}}

\def\absbopb#1{\underline{\underline{#1}}^*}
\def\relbopb[#1]#2{\underline{\underline{#2}}_{|{#1}}^*}
\def\bopb{\@ifnextchar [{\relbopb}{\absbopb}}

\newcommand{\ad}{\operatorname{ad}}

\newcommand{\catMod}{\mathsf{Mod}}
\newcommand{\catFun}[1][]{\mathsf{Hom}_{#1}}

\newcommand{\stkMod}[1][]{\mathfrak{Mod}_{#1}}
\newcommand{\stkSimp}{\mathfrak{S}}
\newcommand{\stkFun}[1][]{\mathfrak{Hom}_{#1}}

\newcommand{\stka}{\mathfrak{A}}

\newcommand{\stkw}{\mathfrak{W}}

\newcommand{\astk}[1]{#1^+}

\newcommand{\OX}{\O_X}

\begin{document}

\title[Deformation-Quantization\ldots]{Deformation-Quantization of Complex Involutive Submanifolds}

\author[A. D'Agnolo]{Andrea D'Agnolo}
\address{Universit{\`a} di Padova\\
Dipartimento di Matematica Pura ed Applicata\\
via G. Belzoni, 7\\ 
35131 Padova, Italy}
\email{dagnolo@math.unipd.it}
\urladdr{www.math.unipd.it/{\~{}}dagnolo/}

\author[P. Polesello]{Pietro Polesello}
\address{Universit{\`a} di Padova\\
Dipartimento di Matematica Pura ed Applicata\\
via G. Belzoni, 7\\ 
35131 Padova, Italy}
\email{pietro@math.unipd.it}

\thanks{The second named author was partially supported by Fondazione Ing.\ Aldo Gini during the preparation of this paper.}

\date{}

\maketitle

\section*{Introduction} 

Let $M$ be a complex manifold, and $T^*M$ its cotangent bundle endowed with the canonical symplectic structure. The sheaf of rings $\shw_M$ of WKB operators provides a deformation-quantization of $T^*M$. On a complex symplectic manifold $X$ there may not exist a sheaf of rings locally isomorphic to $\opb i \shw_M$, for $i\colon X\supset U \to T^*M$ a symplectic local chart. The idea is then to consider the whole family of locally defined sheaves of WKB operators as the deformation-quantization of $X$. To state it precisely, one needs the notion of algebroid stack, introduced by Kontsevich~\cite{Kontsevich2001}. In particular, the stack of WKB modules over $X$ defined in Polesello-Schapira~\cite{Polesello-Schapira} (see also Kashiwara~\cite{Kashiwara1996} for the 
contact case) is better understood as the stack of $\stkw_X$-modules, where $\stkw_X$ denotes the algebroid stack of deformation-quantization of $X$.

Let $V\subset X$ be an involutive (i.e.\ coisotropic) submanifold. 
Assume for simplicity that the quotient of $V$ by its bicharacteristic leaves is isomorphic to a complex symplectic manifold $Z$, and denote by $q\colon V \to Z$ the quotient map. 
If $\shl$ is a simple WKB module along $V$, then the algebra of its endomorphisms is locally isomorphic to $\opb q \opb i \shw_N^\op$, for $i\colon Z\supset U \to T^*N$ a symplectic local chart. Hence we may say that $\shl$ provides a deformation-quantization of $V$. Again, since in general there do not exist globally defined simple WKB modules, the idea is to consider the algebroid stack of locally defined simple WKB modules as the deformation-quantization of $V$. 

In this paper we start by defining what an algebroid stack is, and how it is locally described. We then discuss the algebroid stack of WKB operators on a complex symplectic manifold $X$ and define the deformation-quantization of an involutive submanifold $V\subset X$ by means of simple WKB modules along $V$. Finally, we relate this deformation-quantization  to that given by WKB operators on the quotient of $V$ by its bicharacteristic leaves.

\section{Algebroid stacks}

We start here by recalling the categorical realization of an algebra as in~\cite{Mitchell}, and we then sheafify that construction. We assume that the reader is familiar with the basic notions from the theory of stacks which are, roughly speaking, sheaves of categories. (The classical reference is~\cite{Giraud1971}, and a short presentation is given e.g.\ in~\cite{Kashiwara1996,D'Agnolo-Polesello2003}.)

\medskip

Let $R$ be a commutative ring. An $R$-linear category ($R$-category for short) is a category whose sets of morphisms are endowed with an $R$-module structure, so that composition is bilinear. An $R$-functor is a functor between $R$-categories which is linear at the level of morphisms.

If $A$ is an $R$-algebra, we denote by $\astk A$ the $R$-category with a single object, and with $A$ as set of morphisms. This gives a fully faithful functor from $R$-algebras to $R$-categories.
If $f,g\colon A\to B$ are $R$-algebra morphisms, then transformations $\astk f \Rightarrow \astk g$ correspond to elements $b\in B$ such that $b f(a) = g(a) b$ for any $a\in A$.
Note that the category $\catMod(A)$ of left $A$-modules is $R$-equivalent to the category $\catFun[R](\astk A,\catMod(R))$ of $R$-functors from $\astk A$ to $\catMod(R)$ and that 
the Yoneda embedding
$$
\astk A \to \catFun[R]((\astk A)^\op, \catMod(R)) \approx_R \catMod(A^\op)
$$
identifies $\astk A$ with the full subcategory of right $A$-modules which are free of rank one.
(Here $\approx_R$ denotes $R$-equivalence.)

\smallskip

Let $X$ be a topological space, and $\shr$ a (sheaf of) commutative algebra(s). As for categories, there are natural notions of $\shr$-linear stacks ($\shr$-stacks for short), and of $\shr$-functor between $\shr$-stacks.

If $\sha$ is an $\shr$-algebra, we denote by $\astk\sha$ the $\shr$-stack associated with the separated prestack $U\mapsto \astk{\sha(U)}$. This gives a functor from $\shr$-algebras to $\shr$-categories which is faithful and locally full. If $f,g\colon \sha\to \shb$ are $\shr$-algebra morphisms, transformations $\astk f \Rightarrow \astk g$ are described as above. As above, the stack $\stkMod(\sha)$ of left $\sha$-modules is $\shr$-equivalent to the stack of $\shr$-functors $\stkFun[\shr](\astk \sha,\stkMod(\shr))$, 
and the Yoneda embedding gives a fully faithful functor 
\begin{equation}
\label{eq:Yoneda}
\astk \sha \to \stkFun[\shr]((\astk \sha)^\op,\stkMod(\shr)) \approx_\shr \stkMod(\sha^\op)
\end{equation}
into the stack of right $\sha$-modules. This identifies $\astk\sha$ with the full substack of locally free right $\sha$-modules of rank one.

Recall that one says a stack $\stka$ is non-empty if $\stka(X)$ has at least one object; it is locally non-empty if there exists an open covering $X=\Union\nolimits_i U_i$ such that $\stka|_{U_i}$ is non-empty; and it is
locally connected by isomorphisms if for any open subset $U\subset X$ and any $F,G\in\stka(U)$ there exists an open covering $U=\Union\nolimits_i U_i$ such that $F|_{U_i} \simeq G|_{U_i}$ in $\stka(U_i)$.

\begin{lemma}
\label{le:A+}
Let  $\stka$ be an $\shr$-stack. The following are equivalent
\begin{enumerate}
\item $\stka \approx_\shr \astk \sha$ for an $\shr$-algebra $\sha$,
\item $\stka$ is non-empty and locally connected by isomorphisms.
\end{enumerate}
\end{lemma}

\begin{proof}
By \eqref{eq:Yoneda}, $\astk \sha$ is equivalent to the stack of locally free right
$\sha$-modules of rank one. Then (1) clearly implies (2).
Conversely, let  $\stka$ be an $\shr$-stack as in (2) and $\shp$ an object of 
$\stka(X)$. Then $\sha = \shendo[\stka](\shp)$ is an
$\shr$-algebra and the assignment 
$\shq\mapsto\hom[\stka](\shp,\shq)$ gives an $\shr$-equivalence 
between $\stka$ and $\astk \sha$.
\end{proof}

We are now ready to give a definition of algebroid stack, equivalent to that in Kontsevich~\cite{Kontsevich2001}. It is the linear analogue of the notion of gerbe (groupoid stack) from algebraic geometry~\cite{Giraud1971}.

\begin{definition}
\label{de:stka}
\begin{enumerate}
\item
An $\shr$-algebroid stack is an $\shr$-stack $\stka$ which is locally non-empty and locally connected by isomorphisms.
\item The $\shr$-stack of $\stka$-modules is $\stkMod(\stka) = \stkFun[\shr](\stka,\stkMod(\shr))$.
\end{enumerate}
\end{definition}

Note that $\stkMod(\stka)$ is an example of stack of twisted modules over not necessarily commutative rings (see~\cite{D'Agnolo-Polesello2003}). As above, the Yoneda embedding identifies $\stka$ with the full substack of $\stkMod(\stka^\op)$ consisting of locally free objects of rank one.

\section{Cocycle description of algebroid stacks}

We will explain here how to recover an algebroid stack from local data. The parallel discussion for the case of gerbes can be found for example in~\cite{Breen1994,Breen-Messing2001}.

\medskip

Let $\stka$ be an $\shr$-algebroid stack. By definition, there exists an open covering $U=\Union\nolimits_i U_i$ such that $\stka|_{U_i}$ is non-empty. By Lemma~\ref{le:A+} there are $\shr$-algebras $\sha_i$ on $U_i$ such that $\stka|_{U_i} \approx_\shr \astk{\sha_i}$. Let $\Phi_i \colon \stka|_{U_i} \to \astk{\sha_i}$ and $\Psi_i \colon \astk{\sha_i} \to \stka|_{U_i}$ be quasi-inverse to each other.
On double intersections $U_{ij} = U_i\cap U_j$ there are equivalences $\Phi_{ij} = \Phi_i\Psi_j \colon \astk{\sha_j}|_{U_{ij}} \to \astk{\sha_i}|_{U_{ij}}$.  On triple intersections $U_{ijk}$ there are invertible transformations $\alpha_{ijk} \colon \Phi_{ij}\Phi_{jk} \Rightarrow \Phi_{ik}$ induced by $\Psi_j\Phi_j\Rightarrow\id$. On quadruple intersections $U_{ijkl}$ the following diagram commutes
\begin{equation}
\label{eq:alpha}
\xymatrix{
\Phi_{ij}\Phi_{jk}\Phi_{kl} \ar@{=>}[r]^{\alpha_{ijk}} \ar@{=>}[d]^{\alpha_{jkl}}
& \Phi_{ik}\Phi_{kl} \ar@{=>}[d]^{\alpha_{ikl}} \\
\Phi_{ij}\Phi_{jl} \ar@{=>}[r]^{\alpha_{ijl}} & \Phi_{il} .
}
\end{equation}
These data are enough to reconstruct $\stka$ (up to equivalence), and we will now describe them more explicitly.

On double intersections $U_{ij}$, the $\shr$-functor $\Phi_{ij}\colon \astk{\sha_j} \to \astk{\sha_i}$ is locally induced by $\shr$-algebra isomorphisms. There thus exist an open covering $U_{ij} = \Union\nolimits_\alpha U_{ij}^\alpha$ and isomorphisms of $\shr$-algebras $f_{ij}^\alpha \colon \sha_j \to \sha_i$ on $U_{ij}^\alpha$ such that $\astk{(f_{ij}^\alpha)} = \Phi_{ij}|_{U_{ij}^\alpha}$. 

On triple intersections $U_{ijk}^{\alpha\beta\gamma} = U_{ij}^\alpha \cap U_{ik}^\beta \cap U_{jk}^\gamma$, we have invertible transformations $\alpha_{ijk}|_{U_{ijk}^{\alpha\beta\gamma}} \colon \astk{(f_{ij}^\alpha)}\astk{(f_{jk}^\gamma)} \Rightarrow \astk{(f_{ik}^\beta)}$. There thus exist invertible sections $a_{ijk}^{\alpha\beta\gamma} \in \sha_i^\times(U_{ijk}^{\alpha\beta\gamma})$ such that
$$
f_{ij}^\alpha f_{jk}^\gamma = \ad(a_{ijk}^{\alpha\beta\gamma}) f_{ik}^\beta.
$$

On quadruple intersections $U_{ijkl}^{\alpha\beta\gamma\delta\epsilon\varphi} = U_{ijk}^{\alpha\beta\gamma} \cap U_{ijl}^{\alpha\delta\epsilon} \cap U_{ikl}^{\beta\delta\varphi} \cap U_{jkl}^{\gamma\epsilon\varphi}$,
the diagram \eqref{eq:alpha} corresponds to the equalities
$$
a_{ijk}^{\alpha\beta\gamma} a_{ikl}^{\beta\delta\varphi}
=
f_{ij}^\alpha(a_{jkl}^{\gamma\epsilon\varphi}) a_{ijl}^{\alpha\delta\epsilon}.
$$

Indices of hypercoverings are quite cumbersome, and we will not write them explicitly anymore\footnote{Recall that, on a paracompact space, usual coverings are cofinal among hypercoverings.}.

Let us summarize what we just obtained.

\begin{proposition}
\label{pr:glue}
Up to equivalence, an $\shr$-algebroid stack is given by the following data:
\begin{enumerate}
\item an open covering $X=\Union\nolimits_i U_i$, 
\item $\shr$-algebras $\sha_i$ on $U_i$, 
\item isomorphisms of $\shr$-algebras $f_{ij}\colon \sha_j \to \sha_i$ on $U_{ij}$, 
\item invertible sections $a_{ijk} \in \sha_i^\times(U_{ijk})$, 
\end{enumerate}
such that 
$$
\begin{cases}
f_{ij}f_{jk} = \ad(a_{ijk})f_{ik}, &\text{as morphisms }\sha_k \to \sha_i\text{ on }U_{ijk}, \\
a_{ijk} a_{ikl} = f_{ij}(a_{jkl}) a_{ijl} &\text{in }\sha_i(U_{ijkl}).
\end{cases}
$$
\end{proposition}


\begin{example}
If the $\shr$-algebras $\sha_i$ are commutative, then the 1-cocycle $f_{ij}f_{jk} = f_{ik}$ defines an $\shr$-algebra $\sha$ on $X$ and $\{a_{ijk}\}$ induces a 2-cocycle with values in $\sha^\times$.
In particular, if $\sha_i = \shr|_{U_i}$, then $f_{ij} = \id$ and $\sha=\shr$. Hence, $\shr$-algebroid stacks locally $\shr$-equivalent to $\astk\shr$ are determined by the 2-cocycle $a_{ijk} \in \shr^\times(U_{ijk})$. One checks that two such stacks are (globally) $\shr$-equivalent if and only if the corresponding cocycles give the same cohomology class in $H^2(X;\shr^\times)$.
\end{example}

\begin{example}
Let $X$ be a complex manifold, and $\OX$ its structural sheaf. A line bundle $\shl$ on $X$ is determined (up to isomorphism) by its transition functions $f_{ij}\in\OX^\times(U_{ij})$, where $X = \Union\nolimits_i U_i$  is an open covering such that $\shl|_{U_i} \simeq \O_{U_j}$. Let $\lambda\in\C$, and choose determinations $g_{ij}$ of the multivalued functions $f_{ij}^\lambda$. Since $g_{ij}g_{jk}$ and $g_{ik}$ are both determinations of $f_{ik}^\lambda$, one has $g_{ij}g_{jk} = c_{ijk} g_{ik}$ for $c_{ijk} \in \C_X^\times(U_{ijk})$. 

Let us denote by $\C_{\shl^\lambda}$ the $\C$-algebroid stack associated with the cocycle $\{c_{ijk}\}$ as in the previous example. For $\lambda=m\in\Z$ we have $\C_{\shl^m} \approx_\C \astk{\C_X}$, but in general $\C_{\shl^\lambda}$ is non trivial. On the other hand, $\shl^\lambda$ defines a global object of the algebroid stack\footnote{Here, if $\stka_i$ ($i=1,2$) are $\C$-algebroid stacks, $\stka_1\tens[\C]\stka_2$ is the $\C$-algebroid stack which is equivalent to $\astk {(\sha_1\tens[\C]\sha_2)}$ when $\stka_i \approx_\C \astk \sha_i$, for $\C$-algebras $\sha_i$.} 
$\astk \O_X \tens[\C] \C_{\shl^\lambda}$, so that $\astk \O_X \tens[\C] \C_{\shl^\lambda} \approx_\C \astk \O_X$ is always trivial. Forgetting the $\O$-linear structure, the Yoneda embedding identifies $\shl^\lambda$ with a twisted sheaf in (i.e.\ a global object of) $\stkMod(\C_{\shl^{-\lambda}})$. (Here we used the equivalence $(\C_{\shl^\lambda})^{\op}\approx_{\C}\C_{\shl^{-\lambda}}$.)
\end{example}

\section{Quantization of complex symplectic manifolds}

The relation between Sato's microdifferential operators and WKB operators\footnote{WKB stands for Wentzel-Kramer-Brillouin.} is classical, and is discussed  e.g.~ in~\cite{Pham,AKKT}. We follow here the presentation in~\cite{Polesello-Schapira}, and we refer to~\cite{S-K-K} for the theory of microdifferential operators.

\medskip

Let $M$ be a complex manifold, and denote by $\rho\colon J^1M \to T^*M$ the projection from the 1-jet bundle to the cotangent bundle. 
Let $(t;\tau)$ be the system of homogeneous symplectic coordinates on $T^*\C$, and recall that $J^1 M$ is identified with the affine chart of the projective cotangent bundle $P^*(M\times \C)$ given by 
$\tau\neq 0$.
Denote by $\she_{M\times\C}$ the sheaf of microdifferential operators on $P^*(M\times\C)$.  Its twist by half-forms $\she_{M\times\C}^{\sqrt v} = \opb\pi\Omega_{M\times\C}^{1/2} \tens[\opb\pi\O] \she_{M\times\C} \tens[\opb\pi\O] \opb\pi\Omega_{M\times\C}^{-1/2}$ is endowed with a canonical anti-involution. (Here $\pi\colon P^*(M\times\C)\to M\times\C$ denotes the natural projection.)

In a local coordinate system $(x,t)$ on $M\times\C$, consider the subring $\she_{M\times\C,\hat t}^{\sqrt v}$ of operators commuting with $\partial_t$. 
The ring of WKB operators (twisted by half-forms) is defined by
$$
\shw_M^{\sqrt v} = \oim\rho (\she_{M\times\C,\hat t}^{\sqrt v} |_{J^1M } ).
$$
It is endowed with a canonical anti-involution $*$, and its center is the subfield $k = \shw_{\operatorname{pt}} \subset \C[\![\tau^{-1},\tau]$ of WKB operators over a point.

In a local coordinate system $(x)$ on $M$, with associated symplectic local coordinates 
$(x;u)$ on $T^*M$, a WKB operator $P$ of order $m$
defined on a open subset $U$ of $T^*M$ has a total symbol
$$
\sigma(P)=\sum_{j=-\infty}^m p_j(x;u)\tau^{j},
$$
where the $p_j$'s are holomorphic functions on $U$ subject to the estimates
\begin{equation}\label{eq:estmicrod}
\left\{ \begin{array}{l}
\mbox{for any compact subset $K$ of $U$ there exists a constant}\\
\mbox{$C_K>0$ such that for all $j<0$,}
\sup\limits_{K}\vert p_{j}\vert \leq C_K^{-j}(-j)!.
\end{array}\right.
\end{equation}
If $Q$ is another WKB operator defined on $U$, of total symbol $\sigma(Q)$, then 
$$\sigma(P\circ Q)=\sum_{\alpha\in\N^n} \frac{\tau^{-\vert\alpha\vert}}
{\alpha !} \partial^{\alpha}_u\sigma(P)\partial^{\alpha}_x\sigma(Q).
$$ 

\begin{remark}
The ring $\shw_M^{\sqrt v}$ is a deformation-quantization of $T^*M$ in 
the following sense. Setting $\hbar=\opb \tau$, the sheaf 
of formal WKB operators (obtained by dropping the estimates \eqref{eq:estmicrod}) 
of degree less than or equal to 0 is locally isomorphic to $\O_{T^*M}[\![\hbar]\!]$ as
$\C_{T^*M}$-modules (via the total symbol), and it is equipped with an unitary 
associative product which induces a star-product on $\O_{T^*M}[\![\hbar]\!]$.
\end{remark}

Let $X$ be a symplectic complex manifold of dimension $2n$. Then there are an open covering $X=\Union\nolimits_i U_i$ and symplectic embeddings $\Phi_i \colon U_i \to T^*M$, for $M=\C^n$. Let $\sha_i = \opb{\Phi_i}\shw_M^{\sqrt v}$. Adapting Kashiwara's construction (cf \cite{Kashiwara1996}), Polesello-Schapira~\cite{Polesello-Schapira} proved that there exist isomorphisms of $k$-algebras $f_{ij}$ and invertible sections $a_{ijk}$ as in Proposition~\ref{pr:glue}. Their result may thus be restated as

\begin{theorem}
(cf \cite{Polesello-Schapira})
On any symplectic complex manifold $X$ there exists a canonical $k$-algebroid stack $\stkw_X$ locally equivalent to $\astk{(\opb i \shw_M^{\sqrt v})}$ for any symplectic local chart $i\colon X\supset U \to T^*M$.
\end{theorem}

Note that the canonical anti-involution $*$ on $\shw_M^{\sqrt v}$ extends to an equivalence
of $k$-stacks $\stkw_X\approx_k \stkw_X^\op$. Note also that by Lemma~\ref{le:A+} there exists a deformation-quantization algebra on $X$ if and only if the $k$-algebroid stack $\stkw_{X}$ has a global object, or equivalently if the stack $\stkMod(\stkw_{X})$ has a global object locally isomorphic to
$\opb i \shw_M^{\sqrt v}$ for any symplectic local chart $i\colon X\supset U \to T^*M$.

\section{Quantization of involutive submanifolds}

\medskip

Let $M$ be a complex manifold and $V\subset T^*M$ an involutive\footnote{Recall that $V$ is involutive if for any pair of holomorphic functions $(f,g)$ vanishing on $V$, their Poisson bracket $\{f,g\}$ vanishes on $V$.} submanifold. Similarly to the case of microdifferential operators (for which we refer to \cite{Kashiwara-Oshima} and \cite{Kashiwara1986}), one introduces the sub-sheaf of rings $\shw_V^{\sqrt v}$ of $\shw_M^{\sqrt v}$ generated over $\shw_M^{\sqrt v}(0)$ by the WKB operators $P\in\shw_M^{\sqrt v}(1)$ such that $\sigma_1(P)$ vanishes on $V$. Here $\shw_M^{\sqrt v}(m)$ denotes the sheaf of operators of order less than or equal to $m$, and $\sigma_m(\cdot)\colon \shw_M^{\sqrt v}(m)\to \shw_M^{\sqrt v}(m)/\shw_M^{\sqrt v}(m-1)\simeq \O_{T^*M}\cdot\tau^m$ is the symbol map of order $m$ (which does not depend on the local coordinate system on $M$).

\begin{definition}
\begin{itemize}
     \item[(i)] Let $\shm$ be a coherent $\shw_M^{\sqrt v}$-module. 
     One says that $\shm$ is a regular WKB module along $V$ if locally there exists 
     a coherent sub-$\shw_M^{\sqrt v}(0)$-module $\shm_0$ of $\shm$ which generates it 
     over $\shw_M^{\sqrt v}$, and such that $\shw_V^{\sqrt v}\cdot\shm_0\subset\shm_0$.

     \item[(ii)] One says that $\shm$ is a simple WKB module along $V$ if locally there exists 
     an $\shw_M^{\sqrt v}(0)$-module $\shm_0$ as above such that 
     $\shm_0/\shw_M^{\sqrt v}(-1)\cdot\shm_0\simeq\O_V$.
\end{itemize}
\end{definition}

\begin{example}\label{ex:simple}
Let $(x) = (x_1,\dots,x_n)$ be a local coordinate system on $M$ and
denote by $(x;u) = (x_1,\dots,x_n;u_1,\dots,u_n)$ the associated symplectic local coordinate 
system on $T^*M$. Recall that locally, any involutive submanifold $V$ of codimension $d$ may be
written after a symplectic transformation as:
$$
V = \{(x;u);u_1 = \cdots = u_d = 0\}.
$$
In such a case, any simple WKB module along $V$ is locally isomorphic to
$$
\shw_M^{\sqrt v}/\shw_M^{\sqrt v}\cdot(\partial_{x_1},\dots,\partial_{x_d}).
$$
\end{example}

Let $X$ be a complex symplectic manifold of dimension $2n$, and $V\subset X$ an involutive submanifold. The notions of regular and simple module along $V$ being local, they still make sense in the stack $\stkMod[coh](\stkw_X)$ of coherent WKB modules on $X$. 

\begin{definition}
Denote by $\stkMod[V\text-reg](\stkw_X)$ the full substack  of regular objects along $V$ in $\stkMod[coh](\stkw_X)$, and by $\stkSimp_V$ its full substack of simple objects along $V$.
\end{definition}

Being locally non-empty and locally connected by isomorphisms (see Example \ref{ex:simple}), $\stkSimp_V$ is a $k$-algebroid stack on $X$. Since it is supported by $V$, we consider it as a stack on $V$. 
The first equivalence in the following theorem asserts that the deformation-quantization of $V$ by means of simple WKB modules is equivalent, up to a twist, to that given by WKB operators on the quotient of $V$ by its bicharacteristic leaves.

\begin{theorem}\label{th:simp}
Let $X$ be a complex symplectic manifold, and $V\subset X$ an involutive submanifold. Assume that there exist a complex symplectic manifold $Z$ and a map $q\colon V \to Z$ whose fibers are the bicharacteristic leaves of $V$.
Then there are an equivalence of $k$-algebroid stacks\footnote{Here, if $\stka$ is a $k$-algebroid 
stack on $Z$, $\opb q \stka$ is the $k$-algebroid stack on $V$ which is equivalent to $\astk {(\opb q \sha)}$ when $\stka\approx_k \astk \sha$, for a $k$-algebra $\sha$.} on $V$
\begin{equation}
\label{eq:simp}
\stkSimp_V  \approx_k \opb q \stkw_Z \tens[\C] \C_{\Omega_V^{1/2}},
\end{equation}
and a $k$-equivalence
\begin{equation}
\label{eq:reg}
\stkMod[V\text-reg](\stkw_X) \approx_k \stkMod[coh] ( \opb q \stkw_Z^\op\tens[\C] \C_{\Omega_V^{-1/2}}).
\end{equation}
\end{theorem}

Note that the statement still holds for a general involutive manifold $V$, replacing $\opb q \stkw_Z$ with the algebroid stack obtained by adapting \cite[Proposition 7.3]{Polesello-Schapira} for the symplectic case.

If $V=X$, then $\stkSimp_X\approx_k\stkw_X^\op$ is the stack of locally free left WKB modules of rank one, and $\stkMod[X\text-reg](\stkw_X) \approx_k \stkMod[coh](\stkw_X)$. Since $\Omega_X\simeq\O_X$ by the $n$th power of the symplectic form, one has $\C_{\Omega_X^{1/2}} \approx_\C \astk{\C_X}$. As $q=\id$, the theorem thus reduces to the equivalence $\stkw_X^\op \approx_k \stkw_X$ given by the involution $*$.

If $V = \Lambda$ is Lagrangian, then $Z= \operatorname{pt}$. Hence $\opb q \stkw_{\operatorname pt}  \approx_k  \astk{k_\Lambda}$, and \eqref{eq:simp} asserts that
$$
\stkSimp_\Lambda  \tens[\C] \C_{\Omega_\Lambda^{-1/2}} \approx_k \astk k_\Lambda.
$$
In other words, it asserts that $\stkSimp_\Lambda \tens[\C] \C_{\Omega_\Lambda^{-1/2}} \subset \stkMod(\stkw_X\tens[\C] \C_{\Omega_\Lambda^{1/2}})$ has a global object. This is a result of D'Agnolo-Schapira~\cite{D-S2}, obtained by adapting a similar theorem of Kashiwara~\cite{Kashiwara1996} for microdifferential operators, along the techniques of~\cite{Polesello-Schapira}. As for \eqref{eq:reg}, we recover the WKB analogue of~\cite[Proposition 4]{Kashiwara1996}, also stated in~\cite{D-S2},
$$
\stkMod[\Lambda\text-reg](\stkw_X) \approx_k \stkMod[loc\text-sys] (\astk k_\Lambda\tens[\C]
\C_{\Omega_{\Lambda}^{-1/2}}),
$$
where the right-hand side denotes the stack of twisted locally constant $k$-modules of finite rank.

\begin{proof}[Proof of the theorem]
Consider the two projections
$$
X \from[p_1] X\times Z \to[p_2] Z,
$$
and the associated integral transform functor
$$
\oim{p_1}(\stkMod(\opb{p_2}\stkw_Z^\op) \times \stkMod(\stkw_{X\times Z})) \to \stkMod(\stkw_X).
$$
By the graph embedding, $V$ is identified with a Lagrangian submanifold of $X\times Z$, and $q=p_2|_V$. 
By~\cite{D-S2} there exists a simple module $\shl$ along $V$ in $\stkMod(\stkw_{X\times Z}|_V\tens[\C]\C_{\Omega_V^{1/2}})$. Hence there is an induced $k$-functor
$$
\stkMod(\opb q\stkw_Z^\op\tens[\C]\C_{\Omega_V^{-1/2}}) \to \stkMod(\stkw_X|_V).
$$
This restricts to functors
\begin{align*}
\stkMod[coh] ( \opb q \stkw_Z^\op\tens[\C] \C_{\Omega_V^{-1/2}}) &\to \stkMod[V\text-reg](\stkw_X), \\
\opb q \stkw_Z \tens[\C] \C_{\Omega_V^{1/2}} &\to \stkSimp_V,
\end{align*}
which are local equivalences by the WKB analogue of~\cite[Proposition 4.6]{D-S1}.
\end{proof}

As a corollary, we get a sufficient condition for the existence of a globally defined twisted simple WKB module along $V$.

\begin{corollary}
In the situation of the above theorem, assume that there exists a deformation-quantization algebra on $Z$. Then there exists a globally defined simple object along $V$ in $\stkMod(\stkw_X\tens[\C] \C_{\Omega_V^{1/2}})$.
\end{corollary}

\begin{proof}
Suppose that the $k$-algebroid $\stkw_{Z}$ has a global object. Then its image by the adjunction functor $\stkw_{Z}\to \oim q\opb q \stkw_{Z}$ gives a globally defined twisted simple WKB module along $V$.
\end{proof}

\providecommand{\bysame}{\leavevmode\hbox to3em{\hrulefill}\thinspace}

\end{document}